\documentstyle[12pt,amsfonts]{article}
\input xypic
\newtheorem{theorem}{Theorem}[section]
\newtheorem{corollary}[theorem]{Corollary}
\newtheorem{definition}{Definition}
\newtheorem{example}[theorem]{Example}
\newtheorem{remark}[theorem]{Remark}
\newtheorem{lemma}[theorem]{Lemma}
\newtheorem{proposition}[theorem]{Proposition}
\title{Contra-semicontinuous Functions\thanks{1991 Math.\
Subject Classification --- Primary: 54C08, 54C10, Secondary:
54C05, 54H05. \protect\newline Key words and phrases ---
contra-semicontinuous, semi-regular set, semi-open set, strongly
$S$-closed. \protect\newline Research supported partially by the Japan
- Scandinavia Sasakawa Foundation.\protect\newline This paper was written
while Takashi {\sc Noiri} was visiting the Department of
Mathematics at University of Helsinki in August 1998.}}
\author{Julian Dontchev\\Department of Mathematics\\University
of Helsinki\\PL 4, Yliopistonkatu 15\\00014 Helsinki\\Finland
\and Takashi Noiri\\Department of Mathematics\\Yatsushiro College
of Technology\\2627 Hirayama shinmachi\\Yatsushiro-shi,
Kumamoto-ken\\866 Japan}
\date{}
\begin{document}
\baselineskip=20pt plus 1pt minus 1pt
\newcommand{\fxy}{$f \colon (X,\tau) \rightarrow (Y,\sigma)$}
\newcommand{\ccc}{contra-semicontinuous}
\newcommand{\cccf}{contra-semicontinuous function}
\maketitle
\begin{abstract}
The aim of this paper is to introduce and study the concept of
a \cccf\ and further investigate the class of strongly $S$-closed
spaces. We obtain some new decompositions of generalized
continuous functions.
\end{abstract}

\section{Introduction}\label{s1}

Covering spaces with closed sets has its historical background
in General Topology. 

In 1918, Sierpinski \cite{S1} proved that if a connected compact
Hausdorff space has a countable cover of pairwise disjoint closed
sets, at most one of those sets is nonvoid. In 1992, Cater and
Daily \cite{CD1} showed that if a complete, connected, locally
connected metric space is covered by countably many proper closed
sets, then some two members of these sets must meet. Cater and
Daily improved slightly Sierpinski's result by proving that some
two members must meet in at least continuum many points. Their
new result has applications to several spaces frequently
encountered in functional analysis (see \cite[Corollaries 2,3 and
4]{CD1}). In 1996, the first author \cite{JD1} considered spaces
where every cover by closed sets has a finite subcover. Such
spaces are called strongly $S$-closed. This concept generalizes
Thompson's \cite{T1} $S$-closed spaces, whose definition requires
that every cover by regular closed sets has a finite subcover.
It is a natural to ask which class of generalized continuity
`transforms' strongly $S$-closed spaces onto compact spaces. Such
functions are called in \cite{JD1} contra-continuous. In this
paper we consider a slightly weaker form of contra-continuity
called contra-semicontinuity and study in detail its properties.
We further investigate the class of strongly $S$-closed spaces
and obtain several decompositions of generalized continuous
functions.

\section{Preliminaries}\label{s2}

A subset $A$ of a topological space $(X,\tau)$ is called {\em
semi-regular} \cite{DMN1} if $A$ is both semi-open and
semi-closed. If $A$ is the intersection of an open set and a
semi-closed (resp.\ semi-regular) set, then $A$ is called a {\em
$\cal B$-set} \cite{Tong2} (resp.\ {\em ${\cal A}{\cal B}$-set}
\cite{JD2}). The {\em semi-closure} of $A$, denoted by ${\rm
sCl}(A)$, is the intersection of all semi-closed supersets of
$A$. If ${\rm sCl}(A) \subseteq U$ whenever $U$ is open (resp.\
semi-open), then $A$ is called {\em gs-closed} \cite{AN1} (resp.\
{\em sg-closed} \cite{BL1}). Recall additionally that a set $A$
is called {\em simply-open} \cite{Neu1} if $A = U \cup N$, where
$U$ is open and $N$ is nowhere dense.

The definitions of some basic concepts such as
semi-open set, semi-continuous function, $\alpha$-open set, etc.\
can be found in many papers, for instance \cite{JD1}.

The family of all semi-open (resp.\ semi-closed, semi-regular,
regular open, regular closed, $\alpha$-open, preopen,
$\beta$-open, clopen) subsets of a topological space $(X,\tau)$
will be denoted by $SO(X)$ (resp.\ $SC(X)$, $SR(X)$, $RO(X)$,
$RC(X)$, ${\alpha}(X)$, $PO(X)$, ${\beta}(X)$, $CO(X)$).

\begin{lemma}\label{l1}
For a subset $A$ of a space $X$, the following conditions are
equivalent:

{\rm (1)} $A \in SR(X)$.

{\rm (2)} $A \in {\beta}(X) \cap SC(X)$.
\end{lemma}

\begin{lemma}\label{l2}
For a subset $A$ of a space $(X,\tau)$, the following conditions
are equivalent:

{\rm (1)} $A \in RO(X)$.

{\rm (2)} $A \in \tau \cap SC(X)$.

{\rm (3)} $A \in {\alpha}(X) \cap SC(X)$.

{\rm (4)} $A \in PO(X) \cap SC(X)$.
\end{lemma}

A function \fxy\ is called {\em perfectly continuous} \cite{N1}
(resp.\ {\em completely continuous} \cite{AG1}, {\em
$SR$-continuous}, {\em $RC$-continuous}, {\em $\cal
B$-continuous} \cite{Tong2}, {\em ${\cal A}{\cal B}$-continuous}
\cite{JD2}, {\em simply-continuous} \cite{Neu1})
if the preimage of every open subset of $Y$ is clopen (resp.\
regular open, semi-regular, regular closed, $\cal B$-set, ${\cal
A}{\cal B}$-set, simply-open) in $X$.

A function \fxy\ is called {\em regular set-connected}
\cite{DGR1} (resp.\ {\em $(\theta,s)$-continuous} \cite{JK1},
{\em weakly $\theta$-irresolute} \cite{GNR1}, {\em $R$-map}
\cite{C1}, $\theta$-irresolute \cite{KN1}) if the preimage of
every regular open subset of $Y$ is clopen (resp.\ closed,
semi-closed, regular open, intersection of
regular open sets) in $X$.

\section{Contra-semicontinuous functions}\label{s3}

\begin{definition}\label{d1}
{\em A function \fxy\ is called {\em contra-semicontinuous}
(resp.\ {\em contra-continuous} \cite{JD1}) if the preimage of
every open subset of $Y$ is semi-closed (resp.\ closed) in $X$.}
\end{definition}

The following diagram shows how \ccc\ functions are related to
some similar types of generalized continuity.

$$
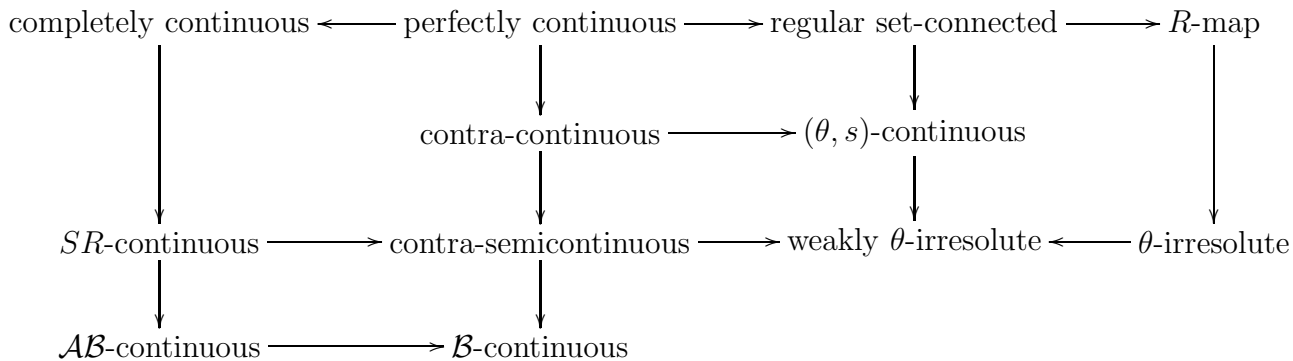
\diagram
\text{completely continuous} \ddto & \text{perfectly
continuous} \lto \dto \rto & \text{regular set-connected}
\dto \rto & \text{$R$-map} \ddto \\ & \text{contra-continuous}
\rto \dto & \text{$(\theta,s)$-continuous} \dto & \\
\text{$SR$-continuous} \rto \dto & \text{contra-semicontinuous}
\dto \rto & \text{weakly $\theta$-irresolute} &
\text{$\theta$-irresolute} \lto \\ \text{${\cal A}{\cal
B}$-continuous} \rto & \text{$\cal B$-continuous}
\enddiagram
$$

\vspace{7mm}

The following examples show that contra-semicontinuity is placed
strictly between contra-continuity and $\cal B$-continuity as
well as strictly between $SR$-continuity and weakly
$\theta$-irresolute continuity.

\begin{example}\label{e1}
{\em A \ccc\ function need not be contra-continuous. Let $f
\colon {\mathbb R} \rightarrow {\mathbb R}$ be the function $f(x)
= [x]$, where $[x]$ is the Gaussian symbol. If $V$ is a closed
subset of the real line, its preimage $U = f^{-1}(V)$ is the
union of intervals of the form $[n,n+1)$, $n \in {\mathbb Z}$;
hence $U$ is semi-open being union of semi-open sets. But $f$ is
not contra-continuous, since $f^{-1}(0.5,1.5) = [1,2)$ is not
closed in $\mathbb R$.}
\end{example}

\begin{example}\label{e2}
{\em The identity function on the real line (with the usual
topology) is $\cal B$-continuous but not contra-semicontinuous,
since the preimage of each singleton fails to be semi-open.}
\end{example}

\begin{example}\label{e3}
{\em A \cccf\ need not be $SR$-continuous. Let $X = \{ a,b,c \}$,
$\tau = \{ \emptyset, \{ a \}, \{ b \}, \{ a,b \}, X \}$ and
$\sigma = \{ \emptyset, \{ c \}, X \}$. The identity function $f
\colon (X,\tau) \rightarrow (X,\sigma)$ is even contra-continuous
but not $SR$-continuous, since $A = \{ a,b \} \in \sigma$ but $A$
is not semi-regular in $(X,\tau)$.}
\end{example}

\begin{example}\label{e4}
{\em A weakly $\theta$-irresolute function need not be \ccc. Let
$X = \{ a,b \}$, $\tau = \{ \emptyset, X \}$ and $\sigma = \{
\emptyset, \{ a \}, X \}$. The identity function $f \colon
(X,\tau) \rightarrow (X,\sigma)$ is weakly $\theta$-irresolute
as only the trivial subsets of $X$ are regular open in
$(X,\sigma)$. However, $f^{-1} (\{ a \}) = \{ a \}$ is not
semi-closed in $(X,\tau)$; hence $f$ is not
contra-semicontinuous.}
\end{example}

\begin{proposition}\label{p1}
For a function \fxy, the following conditions are equivalent:

{\rm (1)} $f$ is \ccc.

{\rm (2)} For every closed subset $F$ of $Y$, $f^{-1} (F) \in
SO(X,\tau)$.

{\rm (3)} For each $x \in X$ and each closed subset $F$ of $Y$
containing $f(x)$, there exists a semi-open $U \in SO(X,\tau)$
such that $f(U) \subseteq F$.

{\rm (4)} ${\rm Int}({\rm Cl} (f^{-1}(V))) = {\rm
Int}(f^{-1}(V))$ for every $V \in \sigma$.

{\rm (5)} ${\rm Cl}({\rm Int}(f^{-1}(F))) = {\rm Cl}(f^{-1}(F))$
for every closed set $F$ of $Y$.
\end{proposition}

Next, we offer the following three decomposition theorems.

\begin{theorem}\label{t1}
For a function \fxy\ the following conditions are equivalent:

{\rm (1)} $f$ is $SR$-continuous.

{\rm (2)} $f$ is $\beta$-continuous and contra-semicontinuous.
\end{theorem}

{\em Proof.} Follows directly from Lemma~\ref{l1}. $\Box$

\begin{corollary}
\cite[Theorem 3.11]{JD1} Every contra-continuous,
$\beta$-continuous is semi-continuous.
\end{corollary}

\begin{example}\label{ee1}
{\em The concepts of $\beta$-continuity and
contra-semicontinuity are independent from each other. Consider
the classical Dirichlet function $f \colon {\mathbb R}
\rightarrow {\mathbb R}$, where $\mathbb R$ is the real line with
the usual topology:

\[ f(x) = \left\{ \begin{array}{ll} 1, &
\mbox{$x \in {\mathbb Q},$} \\ 0, & \mbox{otherwise.}
\end{array} \right. \]

It is easily observed that $f$ is $\beta$-continuous (in fact,
$f$ is even precontinuous). But $f$ is neither contra-continuous
nor $SR$-continuous as $\mathbb Q$ is not semi-closed (hence
not semi-regular).}
\end{example}

\begin{theorem}\label{t2}
For a function \fxy\ the following conditions are equivalent:

{\rm (1)} $f$ is completely continuous.

{\rm (2)} $f$ is precontinuous and contra-semicontinuous.
\end{theorem}

{\em Proof.} Follows directly from Lemma~\ref{l2}. $\Box$

\begin{example}\label{ee2}
{\em The identity function on real line $\mathbb R$ with the
usual topology is (pre)continuous but it is neither \ccc\ nor
completely continuous. For example, $f^{-1} ({\mathbb R}
\setminus \{ 0 \})$ is neither semi-closed nor regular open.}
\end{example}

\begin{theorem}\label{t3}
For a function \fxy\ the following conditions are equivalent:

{\rm (1)} $f$ is $RC$-continuous.

{\rm (2)} $f$ is $\beta$-continuous and contra-continuous.
\end{theorem}

{\em Proof.} Follows easily from the proof of \cite[Theorem
3.10]{JD1}. $\Box$

\begin{remark}
{\em The classical Dirichlet function from Example~\ref{ee1}
shows that a $\beta$-continuous function need not be either
$RC$-continuous nor contra-continuous. On the other hand,
Example~\ref{e3} shows that a contra-continuous function need not
be $\beta$-continuous.}
\end{remark}

\begin{definition}\label{d2}
{\em A function \fxy\ is called {\em contra-gs-continuous}
(resp.\ {\em contra-sg-continuous}) if the preimage of every open
subset of $Y$ is gs-closed (resp.\ sg-closed) in $X$.}
\end{definition}

\begin{theorem}\label{t4}
For a function \fxy\ the following conditions are equivalent:

{\rm (1)} $f$ is \ccc.

{\rm (2)} $f$ is $\cal B$-continuous and contra-gs-continuous.
\end{theorem}

{\em Proof.} (1) $\Rightarrow$ (2) is trivial.

(2) $\Rightarrow$ (1) Let $V \subseteq Y$ be open. Set $f^{-1}(V)
= U \cap F$, where $U \in \tau$ and $F$ is semi-closed in
$(X,\tau)$. Clearly, $f^{-1}(V) \subseteq U$, $U \in \tau$.
Hence, ${\rm sCl} (f^{-1}(V)) \subseteq U$, since $f^{-1}(V)$ is
gs-closed as $f$ is contra-gs-continuous. Now, ${\rm Int}({\rm
Cl}(f^{-1}(V))) = {\rm Int}({\rm Cl}(U \cap F)) \subseteq {\rm
Int}({\rm Cl}(U) \cap {\rm Cl}(F)) = {\rm Int}({\rm Cl}(U)) \cap
{\rm Int}({\rm Cl}(F)) \subseteq {\rm Int}({\rm Cl}(U)) \cap F$,
since $F$ is semi-closed. So, ${\rm Int}({\rm Cl}(f^{-1}(V)))
\cap U \subseteq {\rm Int}({\rm Cl}(U)) \cap U \cap F$. Since
${\rm Int}({\rm Cl}(f^{-1}(V))) \cup f^{-1}(V) = {\rm sCl} 
(f^{-1}(V)) \subseteq U$ and $U \subseteq {\rm Int}({\rm
Cl}(U))$, we have ${\rm Int}({\rm Cl}(f^{-1}(V))) \subseteq U
\cap F = f^{-1}(V)$. This shows that $f^{-1}(V)$ is semi-closed.
$\Box$

The following example will show that the concepts of $\cal
B$-continuity and contra-gs-continuity are independent from each
other and that contra-gs-continuity is strictly weaker than
contra-semicontinuity.

\begin{example}\label{ee3}
{\em Let $X = \{ a,b \}$, $\tau = \{ \emptyset, \{ a \}, X \}$
and $\sigma = \{ \emptyset, X \}$. The identity function on
$(X,\tau)$ is continuous and hence $\cal B$-continuous but not
contra-gs-continuous, because ${\rm sCl} \{ a \} = X
\not\subseteq \{ a \} \in \tau$. Additionally, the identity
function $f \colon (X,\sigma) \rightarrow (X,\tau)$ is
contra-gs-continuous but not $\cal B$-continuous, because
$f^{-1}(\{ a \}) = \{ a \}$ is not a $\cal B$-set in
$(X,\sigma)$.}
\end{example}

The following two results give a tridecomposition of
$SR$-continuity and complete continuity.

\begin{corollary}\label{c1}
For a function \fxy\ the following conditions are equivalent:

{\rm (1)} $f$ is $SR$-continuous.

{\rm (2)} $f$ is $\beta$-continuous, $\cal B$-continuous and
contra-gs-continuous.
\end{corollary}

{\em Proof.} Follows from Theorem~\ref{t1} and Theorem~\ref{t4}.
$\Box$

\begin{corollary}\label{c2}
For a function \fxy\ the following conditions are equivalent:

{\rm (1)} $f$ is completely continuous.

{\rm (2)} $f$ is precontinuous, $\cal B$-continuous and
contra-gs-semicontinuous.
\end{corollary}

{\em Proof.} Follows from Theorem~\ref{t2} and Theorem~\ref{t4}.
$\Box$

\begin{remark}
{\em In both corollaries stated above the three functions in
conditions (2) are pairwise independent. The function $f \colon
(X,\tau) \rightarrow (X,\tau)$ (resp.\ $f \colon (X,\sigma)
\rightarrow (X,\tau)$) in Example~\ref{ee3} is precontinuous but
not contra-gs-continuous (resp.\ not $\cal B$-continuous).}
\end{remark}

\begin{theorem}\label{t5}
For a function \fxy\ the following conditions are equivalent:

{\rm (1)} $f$ is \ccc.

{\rm (2)} $f$ is simply-continuous and contra-sg-continuous.
\end{theorem}

{\em Proof.} It is similar to the proof of Theorem~\ref{t4} and
hence omitted. $\Box$

\begin{example}
{\em Not every simply-continuous function is
contra-sg-continuous. Consider the following function $f \colon
{\mathbb R} \rightarrow {\mathbb R}$, where $\mathbb R$ is the
real line with the usual topology:

\[ f(x) = \left\{ \begin{array}{rl} 1, &
\mbox{$x > 0,$} \\ -1, & \mbox{$x < 0,$} \\ 0, & \mbox{$x=0$.}
\end{array} \right. \]

It can be easily observed that $f$ is simply-continuous. But $f$
is not contra-sg-continuous, since $\{ 0 \}$ is closed and its
preimage $\{ 0 \}$ is not semi-open.}
\end{example}

\begin{example}
{\em Example~\ref{e4} shows that not every contra-sg-continuous
and precontinuous function is simply-continuous. In
Example~\ref{e4}, $f \colon (X,\tau) \rightarrow (X,\sigma)$ is
precontinuous but not contra-sg-continuous.}
\end{example}

\begin{corollary}\label{c3}
For a function \fxy\ the following conditions are equivalent:

{\rm (1)} $f$ is $SR$-continuous.

{\rm (2)} $f$ is $\beta$-continuous, simply-continuous and
contra-sg-continuous.
\end{corollary}

\begin{corollary}\label{c4}
For a function \fxy\ the following conditions are equivalent:

{\rm (1)} $f$ is completely continuous.

{\rm (2)} $f$ is precontinuous, simply-continuous and
contra-sg-semicontinuous.
\end{corollary}

\begin{remark}
{\em The composition of even two contra-continuous functions need
not be \ccc. Let $X = \{ a,b,c \}$, $\tau = \{ \emptyset, \{ a
\}, \{ b \}, \{ a,b \}, X \}$, $\sigma = \{ \emptyset, \{ c \},
X \}$ and $\mu = \{ \emptyset, \{ a,b \}, X \}$. Let $f \colon
(X,\tau) \rightarrow (X,\sigma)$ and $g \colon (X,\sigma)
\rightarrow (X,\mu)$ be the identity functions. Note that both
$f$ and $g$ are contra-continuous but their composition $g \circ
f$ is not even \ccc, since $\{ c \}$ is closed in $\mu$ but
$(g \circ f)^{-1} (\{ c \}) \not\in SO(X,\tau)$.}
\end{remark}

\section{Strongly S-closed spaces}\label{s4}

\begin{definition}
{\em A topological space $(X,\tau)$ is called:

(1) {\em semi-compact} \cite{C1} (resp.\ {\em $s$-closed}
\cite{DMN1}, $S$-closed \cite{T1}) if for every semi-open cover
$\{ V_i \colon i \in I \}$ of $X$, there exists a finite subset
$F$ of $I$ such that $X = \cup \{ V_i \colon i \in F \}$ (resp.\
$X = \cup \{ {\rm sCl} (V_i) \colon i \in F \}$, $X = \cup \{
{\rm Cl} (V_i) \colon i \in F \}$),

(2) {\em nearly compact} \cite{SM1} (resp.\ {\em quasi-H-closed}
\cite{PT1}) if for every open cover $\{ V_i \colon i \in I \}$
of $X$, there exists a finite subset $F$ of $I$ such that $X =
\cup \{ {\rm Int}({\rm Cl}(V_i)) \colon i \in F \}$ (resp.\ $X
= \cup \{ {\rm Cl} (V_i) \colon i \in F \}$),

(3) {\em strongly $S$-closed} \cite{JD1} (resp.\ mildly compact
\cite{St1}) if every closed (resp.\ clopen) cover of $X$ has a
finite subcover.}
\end{definition}

\begin{lemma}
A space $X$ is $s$-closed (resp.\ $S$-closed, nearly compact) if
and only if every semi-regular (resp.\ regular closed, regular
open) cover of $X$ has a finite subcover.
\end{lemma}

The implications in the following diagram are well-known.

$$
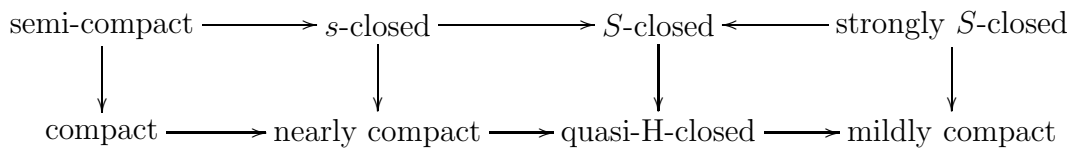
\diagram
\text{semi-compact} \dto \rto & \text{$s$-closed} \dto \rto &
\text{$S$-closed} \dto & \text{strongly $S$-closed} \dto \lto
\\ \text{compact} \rto & \text{nearly compact} \rto & 
\text{quasi-H-closed} \rto & \text{mildly compact}
\enddiagram
$$

\vspace{7mm}

\begin{theorem}\label{t6}
Let \fxy\ be a surjection. If one of the following conditions
holds, then $Y$ is strongly $S$-closed.

{\rm (1)} $f$ is contra-semicontinuous and $X$ is semi-compact,

{\rm (2)} $f$ is $SR$-continuous and $X$ is $s$-closed,

{\rm (3)} $f$ is completely continuous and $X$ is $S$-closed,

{\rm (4)} $f$ is contra-continuous and $X$ is compact,

{\rm (5)} $f$ is $RC$-continuous and $X$ is nearly compact,

{\rm (6)} $f$ is perfectly continuous and $X$ is mildly compact.
\end{theorem}

{\em Proof.} We will prove only the last condition, since proofs
of the other ones are analogous. Let $\{ V_i \colon i \in I \}$
be a closed cover of $Y$. Since $f$ is perfectly continuous, $\{
f^{-1}(V_i) \colon i \in I \}$ is a clopen cover of $X$. Clearly,
there exists a finite $F \subseteq I$ such that $X = \cup_{i \in
F} f^{-1} (V_i)$ as $X$ is mildly compact. Hence, $Y = \cup_{i
\in F} V_i$. This shows that $Y$ is strongly $S$-closed. $\Box$

\section{Some miscellaneous results}\label{s5}

\begin{theorem}\label{51}
Let $(X,\tau)$ be connected and $(Y,\sigma)$ be $T_1$. If \fxy\
is contra-continuous, then $f$ is constant.
\end{theorem}

{\em Proof.} We assume that $Y$ is nonempty. Since $Y$ is a
$T_1$-space, then ${\cal U} = \{ f^{-1} (\{ y \}) \colon y \in
Y \}$ is a disjoint open partition of $X$. If $|{\cal U}| \geq
2$, then there exists a proper nonempty set $W$ (namely, some $U
\in {\cal U}$). Since $X$ is connected, then $|{\cal U}| = 1$.
Hence, $f$ is constant. $\Box$

\begin{corollary}
The only contra-continuous function defined on the real line
$\mathbb R$ is the constant one.
\end{corollary}

Recall that a function \fxy\ is called {\em preclosed}
\cite{DHMN1} if the image of every closed subset of $X$ is
preclosed in $Y$. In 1969, El'kin defined a topological
space $(X,\tau)$ to be {\em globally disconnected} \cite{El1} if
every semi-open set is open. A space $X$ is called {\em locally
indiscrete} if every open set is closed.

\begin{theorem}\label{52}
Let \fxy\ be a contra semi-continuous and pre-closed surjection.
If $X$ is globally disconnected, then $Y$ is locally indiscrete.
\end{theorem}

{\em Proof.} Let $V \in \sigma$. Since $f$ is \ccc, $f^{-1} (V)
= U$ is semi-closed in $X$. Hence $U$ is closed, since $X$ is
globally disconnected. Thus, $f(U) = V$ is preclosed in $Y$ as
$f$ is pre-closed. Now ${\rm Cl}(V) = {\rm Cl}({\rm Int}(V))
\subseteq V$, i.e., $V$ is closed. This shows that $Y$ is locally
indiscrete. $\Box$

\begin{theorem}
Contra-semicontinuos images of hyperconnected spaces are
connected.
\end{theorem}

{\em Proof.} Let \fxy\ be \ccc\ such that $X$ is hyperconnected,
i.e., every open subset of $X$ is dense. Assume that $B$ is a
proper clopen subspace of $Y$. Then $A = f^{-1}(B)$ is both
semi-open and semi-closed as $f$ is \ccc. This shows that $A$ is
semi-regular. Hence, ${\rm Int} (A)$ and ${\rm Int} (X \setminus
A)$ are disjoint nonempty open subsets of $X$. This clearly
contradicts with the fact that $X$ is hyperconnected. Thus, $Y$
is connected. $\Box$

For a topological space $X$, the Cantor-Bendixson derivative
$D(X)$ is the set of all non-isolated points of $X$. A
topological space $(X,\tau)$ is called {\em sporadic} \cite{DGR2}
if the Cantor-Bendixson derivative of $X$ is meager. Recall also
that a space $X$ is called a {\em $T_{\frac{1}{2}}$-space} if
every singleton is open or closed.

\begin{theorem}
If $(X,\tau)$ is CCC (= countable chain condition), $Y$ is a
$T_{\frac{1}{2}}$-space and \fxy\ is \ccc, then $Y$ is sporadic.
\end{theorem}

{\em Proof.} If $Y$ is discrete, then we are done. Assume next
that $D(Y) \not= \emptyset$. Clearly, ${\cal U} = \{ f^{-1} (\{
y \}) \colon y \in D(Y) \}$ is nonempty. Since $f$ is \ccc, $\cal
U$ is a disjoint family of nonempty semi-open subsets of $X$.
Then ${\cal V} = \{ {\rm Int} (f^{-1} (\{ y \})) \colon y \in
D(Y) \}$ is a disjoint family of nonempty open subsets of $X$.
Since $X$ is CCC, $|{\cal V}| \leq \aleph_0$. Hence, $|D(X)| \leq
\aleph_0$. Since every singleton in $D(X)$ is nowhere dense,
$D(X)$ is meager. This shows that $Y$ is sporadic. $\Box$

\
\
E-mail: {\tt dontchev@cc.helsinki.fi}, {\tt
noiri@as.yatsushiro-nct.ac.jp}
\
\
\end{document}